\documentclass[final,3p,10pt]{elsarticle}

\usepackage{graphicx}
\graphicspath{{Figs/}}
\usepackage{placeins}
\usepackage{float}
\usepackage{capt-of}
\usepackage{bm}
\usepackage{algorithm}
\usepackage{algpseudocode}
\usepackage{setspace}
\usepackage{booktabs}
\usepackage{amssymb}
\usepackage{amsmath}
\usepackage{csvsimple}
\usepackage{subcaption}
\usepackage{multirow}

\usepackage[hyphens]{url}
\usepackage[colorlinks=true]{hyperref}
\hypersetup{breaklinks=true}
\newcommand{\review}[1]{\textcolor{black}{#1}}

\usepackage{amsthm}

\usepackage{lineno}

\journal{Elsevier}

\begin{document}


\begin{frontmatter}


\title{Routing Problem for Unmanned Aerial Vehicle \\ Patrolling Missions - A Progressive Hedging Algorithm}



\author[tamu]{Sudarshan Rajan}
\author[lanl]{Kaarthik Sundar}
\author[tamu]{Natarajan Gautam}
\address[tamu]{Texas A\& M University, College Station, TX}
\address[lanl]{Information Systems \& Modeling, Los Alamos National Laboratory, Los Alamos, NM}

\begin{abstract}
The paper presents a two-stage stochastic program to model a routing problem involving an Unmanned Aerial Vehicle (UAV) in the context of patrolling missions. In particular, given a set of targets and a set of supplemental targets corresponding to each target, the first stage decisions involve finding the sequence in which the vehicle has to visit the set of targets. Upon reaching each target, the UAV collects information and if the operator of the UAV deems that the information collected is not of sufficient fidelity, then the UAV has to visit all the supplemental targets corresponding to that target to collect additional information before proceeding to visit the next target. The problem is solved using a progressive hedging algorithm and extensive computational results corroborating the effectiveness of the proposed model and the solution methodology is presented. 
\end{abstract}


\end{frontmatter}


\section{Introduction} \label{sec:introduction}
The use of small Unmanned Aerial Vehicles (UAVs) has seen a tremendous increase in the past decade for both civilian applications \cite{Curry2004,Corrigan2007} viz., precision agriculture \cite{Tokekar2016, Shanahan2001}, forest fire monitoring and management \cite{Casbeer2005}, ocean bathymetry \cite{Ferreira2009} etc. and military applications \cite{Zajkowski2006} viz., monitoring, surveillance, border patrol, intelligence, and reconnaissance (see \cite{Manyam2019,Sundar2017} and references therein). The primary reason for this increase is attributed to the low fixed and operational costs, ease of use, payload capacity, and ability to fly low-altitude missions. There is extensive work in the literature \cite{Shima2009} that examines specific applications and develops algorithms ranging from higher level path planning to lower level control algorithms and the integration of the aforementioned applications. Higher level path planning algorithms are offline algorithms that are used to obtain paths for the vehicles a few hours prior to the start of the mission; in the optimization literature, this corresponds to vehicle routing problems (VRP) for UAVs. The lower level control algorithms are online algorithms that perform real-time trajectory adjustments and vehicle speed changes based on the environment at play. Most of the literature in small UAVs have had an application-specific focus because each application tends to impose its own unique set of constraints that need to be handled separately. One application of small UAVs that has received little attention in literature, in terms of higher level path planning, is that of patrolling \cite{Liu2019}. 

\review{Patrolling applications are essentially data collection missions and despite the term patrolling having a military application connotation, similar application exists even in a civilian context like precision agriculture, crop health monitoring \cite{Tokekar2016}. Such missions are typically associated with uncertainty in the information collected from pre-specified targets. It also entails re-routing to visit additional nearby locations near the pre-specified targets to improve the fidelity of the information gathered. Hence, higher level planning algorithms that take into account this information uncertainty and plan vehicle routes are needed for such applications. There have been several attempts in the last few years to build platforms for operating small UAVs in data gathering missions pertaining to precision agriculture \cite{Tokekar2016}, mapping disaster recovery \cite{Adams2011}, etc. To make things concrete, consider the application of precision farming where the UAV visits a particular target to obtain images of the crops. Suppose the images indicate that the crops at that location are damaged or have an unusual pattern in leaf color, then the mission control would want the vehicle to inspect additional nearby locations to gather more information to come up with robust conclusions on the crop health. If nothing unusual is observed then the vehicle is expected to continue its mission.} In this article, we address the problem of higher level planning that can occur very frequently in almost all data gathering applications that is informally defined as follows: \\

\noindent \textit{Given a set of targets, a source, a destination, and a UAV, the objective is to find a path for the UAV that starts at the source, terminates at the destination, and visits each target exactly once, such that the following conditions are satisfied: if the information (uncertain) collected at a particular target is not sufficient, the UAV has to visit one \cite{rajan2019routing} or more supplemental targets corresponding to the visited target, before proceeding towards the next target in the path and the sum of the cost of the path and the expected additional travel cost to visit the supplemental targets is a minimum. } \\ 

We shall refer to the above problem as the single vehicle data gathering problem (SVDGP). The uncertainty in the SVDGP is associated with the information collected at a particular target. This uncertainty is modeled using a probability distribution and we do not make assumptions on the exact form of the distribution. Rather, we assume that the distribution permits sampling. This makes the problem setting very general and increases its applicability to a very large pool of civilian and military applications detailed in the previous paragraphs. We formulate the SVDGP as a two-stage decision making problem where the first stage (``here-and-now'') decisions are made before the mission under the face of uncertainty and the second stage, recourse decisions are made once the mission starts and uncertainty is revealed. \review{We set up the problem as a two-stage sequential decision making problem instead of a multi-stage Nevertheless for the following reasons: (i) the problem, as stated, has not been addressed previously in the literature, (ii) a solution algorithm to the two-stage version of the problem provides a reasonable indication of whether the algorithm can scale to a multi-stage version, i.e., if the two-stage version of the problem itself is intractable, then the multi-stage version of the problem is definitely not tractable, and finally (iii) a two-stage approximation combined with a rolling horizon algorithm \cite{Balasubramanian2004} could be an efficient way to solve the multi-stage version of the stochastic optimization problem.   
Nevertheless, we note that the formulation and the algorithms presented in this paper can be extended to a multi-stage setting.}

\subsection{Related Work}
The literature contains many variants of UAV routing problems, and algorithms that can obtain optimal solutions and heuristics have been extensively studied for these variants. For ease of exposition, we will analyze the work done in the literature using the following two categories: (1) deterministic UAV routing approaches and (2) stochastic approaches.  

As far as deterministic approaches are concerned, plenty of optimization models and algorithms to compute optimal solutions, and fast heuristics to compute good feasible solutions have received extensive attention over the past decade. As mentioned in the introduction, most of the literature concerning deterministic approaches have had an application-specific focus because each application tends to impose its own unique set of constraints that need to be handled separately. In most of the approaches, the concerned routing problems are modeled as variants of the single traveling salesman problem (TSP), multiple TSP, or Vehicle Routing Problems (VRPs) with additional constraints to model the specific mission at hand; for example see \cite{Sundar2013, Sundar2017JINT, Levy2014}. An interested reader is referred to \cite{Otto2018} for an extensive survey on optimization approaches for deterministic UAV routing problems. To the best of our knowledge, the only work in the literature that addresses the specific application of patrolling is that of \cite{Liu2019} where the authors propose an integer programming formulation to model the patrolling problem as a deterministic multiple TSP and present heuristics to solve the same. 

\review{When compared to the deterministic approaches, literature that focuses on modeling uncertainty in the missions and developing algorithms to solve stochastic variants of the problems concerning UAVs are scarce and very recent. Nevertheless there exists a wide range of literature that deals with stochastic optimization problems in the context of the TSPs and VRPs. An interested reader is referred to \cite{Gendreau1996} for a review of all the variants. For the TSP, stochastic variants involving uncertainty in the targets themselves and uncertainty in travel times are the predominant ones that are dealt with in the literature \cite{Jaillet1988,Carraway1989}. As for the VRPs, stochasticity arising from uncertain demands \cite{Tillman1969,Dror1989}, capacities, customers \cite{Bertsimas1988}, travel times \cite{Dror1989,Jaillet1988a} and various combinations of those \cite{Jezequel1985} have been rigorously examined both from a modeling and algorithm development standpoint. In the context of problems concerning UAVs, two-stage stochastic programming formulations have been explored in single and multiple UAV delivery problems \cite{Shavarani2018,Torabbeigi2018,Venkatachalam2018,Liu2019Stochastic}. To the best of our knowledge, apart from our preliminary conference article that introduces a simpler variant of the SVDGP with only one supplemental location per target \cite{rajan2019routing}, there is no work in the literature that explores such formulations and algorithms in the context of data gathering missions with uncertainty in data collected at points of interest. In \cite{rajan2019routing}, we formulate the SVDGP with only one supplemental location per target to make the second stage optimization problem a linear program; this in turn enables a solution algorithm to compute the optimal two-stage solution using a Benders decomposition. In this paper, we explore a general and more practical setting where there are multiple supplemental locations per target that need to be visited to obtain high fidelity information on the target area under consideration. }

\review{In summary, the contributions of this article are as follows: we develop the first two-stage stochastic programming formulation for a UAV data-gathering mission where uncertainty arises in the actual data collected in the targets. The formulation, as presented in the subsequent sections, contains binary variables in both the first as well as the second stage and we present a progressive hedging algorithm \cite{Watson2011} to solve this problem. The motivation behind developing a progressive hedging algorithm to solve this two-stage problem formulation will be detailed in the later sections. Finally, extensive computational experiments that corroborate the effectiveness of these models against deterministic counterparts and effectiveness of the algorithm to compute an optimal solution for any instance of the problem, as opposed to solving an extensive form (EF) using sample average approximation (SAA), are presented. } The rest of the article is organized as follows: in Sec. \ref{sec:pr-statement}, we present the formal problem statement after introducing the necessary notations. In Sec. \ref{sec:formulation}, we formulate the SVDGP as a two-stage stochastic program and present an algorithm to solve the same in Sec. \ref{sec:solution}. Finally in Sec. \ref{sec:results}, we present the computational results followed by conclusions and future research directions in Sec. \ref{sec:conclusion}.

\section{Problem Statement} \label{sec:pr-statement}
Before we present the formal problem statement, we introduce some notations that will be used throughout the rest of the article. We are given a set of $n$ targets $\hat T = \{t_1, \dots, t_n\}$, a source $s$, and a destination $d$. For ease for exposition, we will assume that $s = d = t_0$ and remark that the formulations and algorithms can be extended easily to the case where they are distinct. We also let $T = \hat T \cup \{t_0\}$. As detailed in the previous section, the decision making process is two-staged. The first stage decision for the mission is to compute a path for the UAV that starts and ends at $t_0$ and visits each target in the set $\hat T$ and collects information. Also associated with each target $k \in T$ is a set of supplemental locations $S_k$ whose purpose is as follows: suppose the UAV visits the target $k$ and collects some information, then depending on the this information gathered, it may have to visit some additional set of locations (supplemental location set $S_k$) to either validate the collected information or to gather more information. The decision on whether the UAV needs to visit these supplemental targets are recourse decisions and are taken after the uncertainty is realized. For the purpose of this article, we assume $|S_k|$, the cardinality of $S_k$ for each $k=1, \dots, n$, is $m$, where $m \geqslant 2$ and $S_0 = \emptyset$. We shall, from here on, refer to the set $V = T \bigcup S_1 \bigcup \dots \bigcup S_n$ as the set of vertices. Given these notations, the SVDGP is formulated on a graph $G = (V, E)$ where $E$ denotes the edge set; we delegate the definition of the edge set $E$ to the Sec. \ref{sec:formulation}. We note that $G$ can be a directed or an undirected graph and the formulation and the algorithms presented in this article are applicable to both cases. To keep the exposition fairly general, we will assume that the graph is directed and that the cost of traveling between any pair of targets is asymmetric i.e., given $(i, j) \in V$, $c_{ij} \neq c_{ji}$. The cost of traversing an edge $(i, j) \in E$ can be anything ranging from distance to travel time. We also assume that this cost of traversal between any pair of vertices is known a-priori or can be pre-computed. Since the focus of the article is to develop a model to account for uncertainty in a systematic way, in the next section, we detail the uncertainty model associated with the information collected at each target. 

\subsection{Uncertainty Modeling} \label{subsec:uncertainty}
The uncertainty in the problem arises from the fidelity of information collected at any target $i \in T$. Based on the information collected at target $i$, the UAV may or may not have to visit the supplemental locations in the set $S_i$ before proceeding to visit the next target in the route. Hence, we associate with each target $i \in T$, a Bernoulli random variable with probability $p_i$ which denotes the probability that the information collected at target $i$ is not of sufficient fidelity. For the target $t_0$, we assume that $p_i = 0$, where $i \in \{t_0\}$. For ease of exposition, we assume that the random variables corresponding to any pair of targets, is independent of each other. However, this can be relaxed or changed and the uncertainty can also be modeled using a Markov chain. The uncertainty in this information collected over all the targets is modeled as a binary scenario vector of size $|T|$ where the component corresponding to target $i \in T$ takes a value $1$ with a probability $p_i$. These realizations of uncertainty are captured in a countable set of scenarios $\Omega$, in which each element $\omega \in \Omega$ occurs with probability $p_{\omega}$. Each element $\omega \in \Omega$ is a $|T|$-dimensional binary vector where each element is $1$ with a probability $p_i$. We let $\omega_i$ denote the $i$\textsuperscript{th} component of $\omega$. \review{When $\omega_i = 1$, the data gathered from the target $i$ is not of sufficient fidelity and additional supplemental targets corresponding to the target $i$ have to be visited by the vehicle.} In this particular case, given that the probabilities $p_i$ and $p_j$ associated with any pair of targets $(i,j)$ are independent, the probability of occurrence of the scenario $\omega$ is given by
\begin{flalign}
p_{\omega} = \prod_{i=0}^{n} \left\{ \mathbf{1}(\omega_i=1)\cdot p_i + \mathbf{1}(\omega_i=0) \cdot (1-p_i)\right\} \label{eq:prob_omega}
\end{flalign}
where, the function $\mathbf{1}(\cdot)$ is the indicator function. 

If we let $\mathcal F$ denote the set of all feasible first stage paths for the UAV, the first stage cost is the cost of the path and for any feasible path in $\mathcal F$, the second stage decision is a set of recourse decisions which is a set of routes through the supplemental targets based on the fidelity of information collected at the targets. The second stage decisions are scenario-dependent and the second stage cost is given by the expected additional cost of traversal for the UAV to visit the supplemental targets. The goal of the SVDGP is to find a path that minimizes the sum of the first and second stage travel costs. 

\section{Mathematical Formulation} \label{sec:formulation}
We remark that when there is no uncertainty associated with the problem and when the information gathered at every target $i \in T$ is already of sufficient fidelity, then no supplemental targets need to be visited by the vehicle and the first stage problem reduces to a asymmetric TSP on the set $T$. For the SVDGP, if the information associated with a target $i \in T$ is not of sufficient fidelity, the UAV needs to visits the supplemental targets in the set $S_i$ before proceeding towards the next target. To formulate this problem as a two-stage stochastic program, we first define two edge sets $E^1$ and $E^2$ for the first and the second stages respectively. The edge set $E^1$ includes edges that are allowed for the vehicle in the first stage i.e., any edge between every pair of vertices in the set $T$. The second stage edge set $E^2$ includes the edges that the vehicle is allowed to traverse in the second stage i.e., apart from containing all the edges in the set $E^1$, it includes the edges from every target $i \in T$ to every supplemental target in the set $S_i$, from every supplemental target in the set $\cup_i S_i$ to every target in the set $T$ and for every $i \in T$, between any pair of supplemental targets in the set $S_i$. For any edge in the set $(i, j) \in E^1 \cup E^2$, $c_{ij}$ denotes the cost of traversal of that edge. Also, given a subset of vertices $\hat V$, we define $\delta^1_+(\hat V) = \{(i, j) \in E^1: j \notin \hat{V} \text{ and } i \in \hat V\}$ and $\delta^2_+(\hat V) = \{(i, j) \in E^2: j \notin \hat V \text{ and } i \in \hat V\}$ as the set of outgoing edges in the first and second stage edge set, respectively. Similarly, given $\hat V \subset V$, we let $\delta^1_-(\hat V) = \{(j, i) \in E^1: j \notin \hat V \text{ and } i \in \hat V\}$ and $\delta^2_-(\hat V) = \{(j, i) \in E^2: j \notin \hat V \text{ and } i \in \hat V\}$ as the set of incoming edge. Furthermore, when $\hat V = \{i\}$ i.e., a singleton, we simply write $\delta^1_+(i)$ instead of to $\delta^1_+(\{i\})$.  Finally, given two disjoint subsets of vertices $V_1, V_2 \subset V$, we define $\gamma(V_1 \rightarrow V_2) = \{(i, j) \in E^2: i \in V_1, j \in V_2\}$. 

\subsection{Objective function} \label{subsec:obj}
Given the above notations, we introduce binary first stage decision variables $x_{ij}$ for each $(i,j) \in E^1$, denoting the presence of the edge $(i,j)$ in the first stage solution and another set of binary second stage decision variables $y_{ij}^{\omega}$ for each $(i, j) \in E^2$ denoting the presence of the edge $(i,j)$ in the second stage for the scenario $\omega \in \Omega$. We let $\bm x$ and $\bm y$ denote the vector of first and second stage decision variables, respectively. Finally, given a subset of edges $\hat E^1 \subset E^1$ and $\hat E^2 \subset E^2$, we let $x(\hat{E}^1)$ and $y^{\omega}(\hat{E}^2)$ denote the sums $\sum_{(i,j) \in \hat E^1} x_{ij}$ and $\sum_{(i,j) \in \hat E^2} y_{ij}^{\omega}$, respectively. Now, the objective for the two-stage stochastic programming formulation of the SVDGP is given by:

\begin{flalign}
\min~ C & \triangleq \sum_{(i, j) \in E^1} c_{ij} x_{ij} + \mathbb E_{\Omega} \left[ \beta(\bm x, \omega) \right] \label{eq:obj-full1} \\
& = \sum_{(i, j) \in E^1} c_{ij} x_{ij} + \sum_{\omega \in \Omega} p_{\omega} \beta(\bm x, \omega) \nonumber
\end{flalign} 
Here, $\beta(\bm x, \omega)$ is the additional traversal cost required to visit the supplemental targets given a scenario $\omega \in \Omega$. 

\subsection{First stage constraints} \label{subsec:first stage-constraints} 
Since the first stage solution is a feasible tour through the set of targets in $T$, the constraints in this stage are the TSP constraints which are given below:
\begin{subequations}
\begin{flalign}
x(\delta^1_+(i)) = 1 &\quad \forall i \in T, \label{eq:outdegree-1} \\ 
x(\delta^1_-(i)) = 1 &\quad \forall i \in T, \label{eq:indegree-1} \\ 
x(\delta^1_+(S)) \geqslant 1 &\quad \forall S \subset T, ~ |S| \geqslant 2, \text{ and } \label{eq:sec-1} \\ 
x_{ij} \in \{0, 1\} &\quad \forall (i, j) \in E^1. \label{eq:binary-1}
\end{flalign}
\label{eq:first stage}
\end{subequations}
\noindent Here, Eqs. \eqref{eq:outdegree-1} and \eqref{eq:indegree-1} are the in-degree and out-degree constraints that enforce exactly one edge to enter and leave each target in the set $T$. Eq. \eqref{eq:sec-1} ensures that there are no sub-tours in the solution, and finally Eq. \eqref{eq:binary-1} enforces the binary restrictions on the decision variables $x_{ij}$. 

\subsection{Second stage formulation} \label{subsec:second stage}
The second stage model for a fixed first stage solution $\bm x$ and the realization of uncertainty $\omega \in \Omega$ is as follows:
\begin{subequations}
\begin{flalign}
\beta(\bm x, \omega) = \min \quad \sum_{(i, j) \in E^2} c_{ij} y_{ij}^{\omega} - &\sum_{(i, j) \in E^1} c_{ij} y_{ij}^{\omega}  \label{eq:obj-2}  \\
\quad \text{ subject to:} \qquad \qquad \qquad & \notag \\
y^{\omega}(\delta^2_+(i)) = 1 &\quad \forall i \in T, \label{eq:indegree-2}\\ 
y^{\omega}(\delta^2_-(i)) = 1 &\quad \forall i \in T, \label{eq:outdegree-2} \\ 
y^{\omega}(\delta^2_+(j)) = \omega_i &\quad \forall j \in S_i, i \in T, \label{eq:indegree-supplemental-2}\\ 
y^{\omega}(\delta^2_-(j)) = \omega_i &\quad \forall j \in S_i, i \in T \label{eq:outdegree-supplemental-2}\\
y_{ij}^{\omega} = x_{ij}(1-\omega_i) &\quad \forall (i,j) \in E^1 \label{eq:on-off-1}\\ 
y^{\omega}(\gamma(S_i \rightarrow \{j\})) = x_{ij} \omega_i &\quad \forall i, j \in T, \label{eq:on-off-2} \\
y^{\omega}(\delta^2_+(S)) \geqslant 1 &\quad \forall S \subset V, ~ |S| \geqslant 2, S \cap T \neq \emptyset \text{ and } \label{eq:sec-2}\\
y_{ij}^{\omega} \in \{0, 1\} &\quad \forall (i, j) \in E^2 \label{eq:binary-2}.
\end{flalign}
\label{eq:second stage}
\end{subequations}
The objective in Eq. \eqref{eq:obj-2} minimizes the additional cost of visiting the supplemental targets if information collected at a particular target is not of sufficient fidelity. The first summation is the total cost of the route including the additional supplemental target visits and the second term is the cost of the edges that are present both in the first and the second stage solutions. The Eqs. \eqref{eq:indegree-2} and \eqref{eq:outdegree-2} are the in-degree and out-degree constraints for the targets in the set $T$. The Eqs. \eqref{eq:indegree-supplemental-2} and \eqref{eq:outdegree-supplemental-2} ensure that if the information collected at any target $i \in T$ is not of sufficient fidelity i.e., $\omega_i = 1$, then the supplemental targets in the set $S_i$ have to be visited by the vehicle. The Eq. \eqref{eq:on-off-1} ensures that if an edge $(i,j) \in E^1$ is traversed by the vehicle in the first stage and if $\omega_i = 0$, then this edge is necessarily used in the second stage by setting $y_{ij}^{\omega} = 1$ and if an edge $(i,j) \in E^1$ is not used by the vehicle in the first stage, then it is also not used in the second stage corresponding to scenario $\omega$. Eq. \eqref{eq:on-off-2} ensures that given a target $i\in T$ with $\omega_i = 1$ and given that $j \in T$ is the target that vehicle visits immediately after $i$ in the first stage solution, the vehicle has to collect additional information from all the supplemental targets in the set $S_i$ before visiting the next target $j$. Together, Eqs. \eqref{eq:on-off-1} and \eqref{eq:on-off-2} also ensure that the sequence in which the vehicle visits the targets in the first stage solution given by $\bm x$ remains unchanged in the second stage. The Eq. \eqref{eq:sec-2} eliminate sub-tours in the second stage solution and finally, the Eq. \eqref{eq:binary-2} impose binary restrictions on the decision variables $y_{ij}^{\omega}$. 

\section{Solution Methodology} \label{sec:solution} 
This section details the description of our approach to solve the two-stage mathematical formulation of the SVDGP presented in Sec. \ref{sec:formulation}. The formulation, as presented, contains binary decision variables both in the first and the second stages. We leverage the Progressive hedging (PH) algorithm proposed by Rockafellar and Wets \cite{Rockafellar1991} to solve the SVDGP; PH is also referred to as scenario decomposition in the literature \cite{Watson2011} since it decomposes stochastic programs by scenarios i.e., samples of the realization of the random variables in the problem. \review{The other predominantly used method for solving mixed-integer stochastic optimization problems is that of Stochastic Dual Dynamic Programming (SDDP) \cite{Shapiro2011}. This method relies on outer-approximating the recourse problems using linear constraints obtained from the dual values of linear relaxation of the recourse problems. For the SVDGP, we favour the PH algorithm over the SDDP since each recourse problems itself has exponential number of sub-tour elimination constraints that are added dynamically during the solve and obtaining dual values of the dynamically added constraints requires separate cut-management algorithm. PH is a well-known algorithm for solving multistage stochastic convex optimization problems \cite{Atakan2018}. In fact, PH possesses theoretical convergence guarantees when all decision variables in the convex multistage stochastic program are continuous, i.e., PH is guaranteed to converge to the global optimal solution of the multistage stochastic program when all but the here-and-now decision variable are continuous. In the presence of discrete variables, a wealth of recent theoretical and empirical research \cite{Fan2010,Listes2005,Lokketangen1996} has shown that the PH algorithm can prove to be a very robust heuristic to solve stochastic programs, specifically the case of pure binary programs in both stages of the formulation even in the case that the number of scenarios is prohibitively large. These algorithmic features of PH make its application ideal for the the SVDGP.} To that end, the subsequent sections detail the techniques involved in the PH approach and present findings that shows its impact as an effective heuristic to solve the SVDGP. In the forthcoming paragraphs, we present an overview of the PH algorithm specifically tailored for the SVDGP. 

\subsection{Algorithm Overview} \label{subsec:ph}
Before, we present the overview of the algorithm, we restate the two-stage stochastic program for the SVDGP in a concise manner as follows:
\begin{subequations}
\begin{flalign}
\min \quad & \sum_{(i, j) \in E^1} c_{ij} x_{ij} + \sum_{\omega \in \Omega} p_{\omega} \beta(\bm x, \omega) \label{eq:obj-concise} \\ 
\text{subject to: } \quad & \bm x \in \mathcal Q \label{eq:constraint-set}
\end{flalign}
\label{eq:concise}
\end{subequations}
where, $\mathcal Q = \{ \bm x: \bm x \text{ satisfies Eq. \eqref{eq:first stage}} \}$. The concise formulation, as presented in Eq. \eqref{eq:concise}, is the well-known extensive form of the two-stage stochastic program  \cite{Wallace2005} in Sec. \ref{sec:formulation}. The PH is a scenario-decomposition algorithm that uses a separate set of first stage decision variables for each scenario to perform parallel solves i.e., for each $\omega \in \Omega$, it introduces decision variables $\bm x_{\omega}$ and implicitly enforces the non-anticipativity constraints ($\bm x = \bm x_{\omega} ~~\forall \omega \in \Omega$) via penalization; these constraints avoid allowing the first stage decision vector $\bm x$ to depend on the scenario. The basic PH algorithm takes as input two parameters (i) a penalty factor, $\rho > 0$ and (ii) a termination threshold, $\varepsilon$. Given $\rho$ and $\varepsilon$, the pseudo-code for the PH algorithm is as follows:

\begin{algorithm}
    \caption{Progressive Hedging: A pseudo-code}
    \label{algo:ph}{
    \onehalfspacing
    \begin{algorithmic}[1] 
        \Statex \textbf{Initialization:} 
        \State $k \gets 0$ \Comment{iteration count}
        \State \label{step:main0} For every $\omega \in \Omega$, $\bm x_{\omega}^k \gets \operatornamewithlimits{argmin}_{\bm x \in \mathcal Q} \sum_{(i, j) \in E^1} c_{ij} x_{ij} + \beta(\bm x, \omega)$ 
        \State $\bm x^k \gets \sum_{\omega \in \Omega} p_{\omega} \bm x_{\omega}^k$ 
        \State For every $\omega \in \Omega$, $\bm w_{\omega}^k \triangleq \rho \cdot (\bm x_{\omega}^k - \bm x^k)$ \Comment{initial weight vector computation $\bm w_{\omega}^k$}
        \Statex \textbf{Iteration update:}
        \State \label{step:goto} $k \gets k+1$ 
        \Statex \textbf{Decomposition:}
        \State \label{step:main} For every $\omega \in \Omega$, $$\bm x_{\omega}^k \gets \operatornamewithlimits{argmin}_{\bm x \in \mathcal Q} \left( \sum_{(i, j) \in E^1} c_{ij} x_{ij} + \langle\bm w_{\omega}^{k-1}, \bm x \rangle + \frac{\rho}2 \left\| \bm x - \bm x^{k-1} \right\|^2 + \beta(\bm x, \omega) \right)$$ 
        \Statex \textbf{Aggregation:}
        \State $\bm x^k \gets \sum_{\omega \in \Omega} p_{\omega} \bm x_{\omega}^k$ 
        \Statex \textbf{Weight update:}
        \State For every $\omega \in \Omega$, $\bm w_{\omega}^k \triangleq \bm w_{\omega}^{k-1} + \rho \cdot (\bm x_{\omega}^k - \bm x_{\omega})$
        \Statex \textbf{Termination criterion check:}
        \State $\epsilon^k \triangleq \sum_{\omega \in \Omega} p_{\omega} \left\| \bm x_{\omega}^k - \bm x^k \right\|$ 
        \If{$\epsilon^k > \varepsilon$} 
            \State Go to Step \ref{step:goto}
        \Else
            \State Terminate with $\bm x^k$ as the first stage solution
        \EndIf
    \end{algorithmic}
    }
\end{algorithm}

In Step \ref{step:main} of Algorithm \ref{algo:ph}, the symbol $\langle a, b \rangle$ denotes the dot product of the vectors $a$ and $b$. Furthermore, we notice that the computations in Steps \ref{step:main0} and \ref{step:main} involve solving multiple mixed-integer linear programs, one for each scenario in the set $\Omega$ and that the solves are completely parallelizable. The quadratic penalty term (or the proximal term \cite{Watson2011}) penalizes the non-anticipativity constraints and ensures that it is satisfied as the algorithm terminates. \review{For every iteration $k$, the term $\epsilon^k$ is also referred to as the residual. The values of of $\varepsilon$ is a tolerance value and it can be set to any value that is desired depending on the level of accuracy that is required for PH algorithm. As for the value of $\rho$, the following procedure was used to set its value. This procedure is not new and has been previously used in the literature concerning PH algorithms in \cite{Watson2011}. It was observed in \cite{Watson2011} that values of $\rho \in (0, 1)$ is a good starting point to fix the value of $\rho$ and that $\rho$ values proportional to the first stage edge costs, i.e., $c_{ij}$ would lead to good convergence behaviour (without oscillations) of the PH algorithm. This technique of fixing a $\rho$ value proportional to the costs is referred to as ``cost proportional'' heuristic \cite{Watson2011}. For the SVDGP, we combined both the suggestions by first scaling all the costs between every pair of targets in the range $(0, 1)$ and then allowing the value of $\rho$ to be the average of all the scaled cost values. This technique worked well for all the computational experiments and hence, no further experiments were done to compute more effective values of $\rho$.}

\review{In general, there is no guarantee that the Progressive Hedging algorithm will converge to a binary first stage solution \cite{Rockafellar1991}. But, it has been observed that for a wide variety of problems it is observed to converge to a binary first stage solution \cite{Watson2011}. In all our computational experiments, we did not find a single instance for which the algorithm did not converge. Nevertheless, heuristics like variable-fixing, slamming, or rounding can be used to force convergence.} The following sections show the results obtained by adopting this approach to our two-stage formulation of the problem.

\section{Computational Results} \label{sec:results} 
The PH algorithm for the SVDGP was implemented  using the C++ programming language, and CPLEX 12.8 was used as the underlying solver to solve the multiple MILPs that occur in each iteration of the PH algorithm. All the computational experiments were performed using an 2.9 GHz Intel Core i7 processor with 16 GB RAM. The performance of the algorithm was tested on randomly generated test instances. The instances generation procedure is detailed in the subsequent section. 

\subsection{Instance Generation} \label{subsec:instance-generation} 
All the targets, the supplemental locations, the source, and the destination were all generated on a $100 \times 100$ grid. The source and the destination vertices for all the instances were located at $(5, 5)$ and $(95, 95)$, respectively. The number of targets $n$ were varied from $10$ to $40$ in steps of $5$. For each target, the locations of its corresponding supplemental targets were also randomly generated within a maximum pre-specified radius, $R$, from the target location; the value of $R$ was chosen from the set $\{5, 10\}$ units. The number of supplemental locations per target, $m$ was chosen from the set $\{3, 5, 7, 9\}$. Now, as for the vehicle itself, we assume that the vehicle is a fixed-wing UAV with a minimum turn-radius of $5$ units. For this fixed-wing UAV, the cost of traversal between any pair of targets is assumed to be the shortest path taken by the vehicle to go from one target to the other. This path length is in-turn a function of the heading angle of the UAV at each target. To that end, we generate random heading angles for every target and its corresponding supplemental location and then compute the length using the well-known result by Dubins \cite{Dubins1957}. We remark that the shortest path computed using the result in \cite{Dubins1957} is asymmetric. Given this instance generation procedure, the total number of instances was $56$ and all computational experiments were performed on this set of $56$ instances. A subset of these instances are used for each of the computational experiments and the subset of chosen instances for each experiment is presented in the respective sections. 

\subsection{PH Algorithm Parameters and Scenario Generation} \label{subsec:scenario-generation}
The PH algorithm has two main input parameters, $\rho$ and $\epsilon$. The value of $\epsilon$ is set to $1 \times 10^{-5}$ and the value of $\rho$ for each instance was computed using the heuristic detailed in Sec. \ref{subsec:ph}. We note that for majority of the instances the value of $\rho$ computed using the heuristic was close to $0.5$. The uncertain scenarios for each run of the problem are generated using a vector of Bernoulli random variable, one for each target. Across a majority of the experiments in the subsequent sections, the probability that the information collected at a particular target is of sufficient fidelity is assumed to be $0.5$ for every target. For a few other computation experiments other probability values are used and these values are presented when the corresponding experiment is described. For the source and the destination, we set this probability to $0$ and required number of scenarios are generated from this vector of Bernoulli random variables. The number of scenarios was varied from a minimum value of $10$ to a maximum value of $200$ across all experiments. Each randomly generated scenario, $\omega$ is a binary vector indicating if the information collected at a particular target is of sufficient fidelity or not. In the subsequent paragraphs, we present the results for the various computational experiments perform to evaluate the PH algorithm on the SVDGP. Before we present the results obtained using the PH algorithm, we remark that when the full two-stage stochastic formulation in Sec. \ref{sec:formulation} was provided to CPLEX with a time-limit of three hours, all problem instances with number of targets greater than $50$ and number of scenarios greater than $20$ timed-out. Hence, we do not present any results that show that CPLEX was not able to solve the full problem as stated in Sec. \ref{sec:formulation}. 

\review{\paragraph{Performance of PH algorithm} In this section, we compare the performance of the PH algorithm against solving the full two-stage problem in its extensive form using CPLEX. For this study,  we choose only the 10 target instances. The values of $m$ and $R$ are varied in the set $m \in \{3, 5, 7, 9\}$ and $R \in \{5, 10\}$. Also, the number of scenarios, $|\Omega|$ was varied in the set $\{50, 100, 200\}$ and the probabilities to generate these scenarios were chosen in the set $\{0.2, 0.5, 0.8\}$. In total, the number of test instances for this study was $720$. CPLEX was used to solve the extensive form of the two-stage problem with a computation time limit of $7200$ seconds. The Table \ref{tab:phvscplex} shows the number of instances where solving the extensive form failed to result in an optimal solution for the instance compared to solving the problem using the PH algorithm.  As seen from Table \ref{tab:phvscplex}, it is clear that solving the extensive form using CPLEX is very inefficient since many instances with just $10$ targets were not solved to optimality. On the other hand, the PH algorithm is capable of solving all instances to optimality. This justifies the need for the proposed decomposition algorithm for the SVDGP. For all the instances that were solved to optimality by both the solution methods, the difference in solution quality between the two methods was found to be $0\%$. This result provides an empirical proof that despite the PH algorithm being an heuristic approach, it is capable of producing the globally optimal solution for the two-stage problem. Finally, for all the runs where both the approaches provided the optimal solution, the PH algorithm was on an average $2.5$ times faster than solving the extensive form using CPLEX. Furthermore, out of the $720$ instances, $48\%$ of the instances had PH algorithm converging to the optimal solution in more than $1$ iteration.}

\begin{table}[htbp]
    \centering{
    \onehalfspacing
    \begin{tabular}{cccc}
        \toprule 
        \multirow{2}{*}{$|\Omega|$} & \multirow{2}{*}{total \# instances} & \multicolumn{2}{c}{\# timed out instances} \\
        \cmidrule{3-4}
        & & extensive form & PH algorithm \\
        \midrule 
        50 & 240 & 21 & 0\\
        100 & 240 & 18 & 0\\
        200 & 240 & 42 & 0\\
        \bottomrule
    \end{tabular}
    \caption{Number of instances where providing the extensive form directly to CPLEX timed out as against applying the PH algorithm to solve the two-stage problem.}
    \label{tab:phvscplex}
    }
\end{table}

\begin{figure}
    \centering
    \includegraphics[scale=0.5]{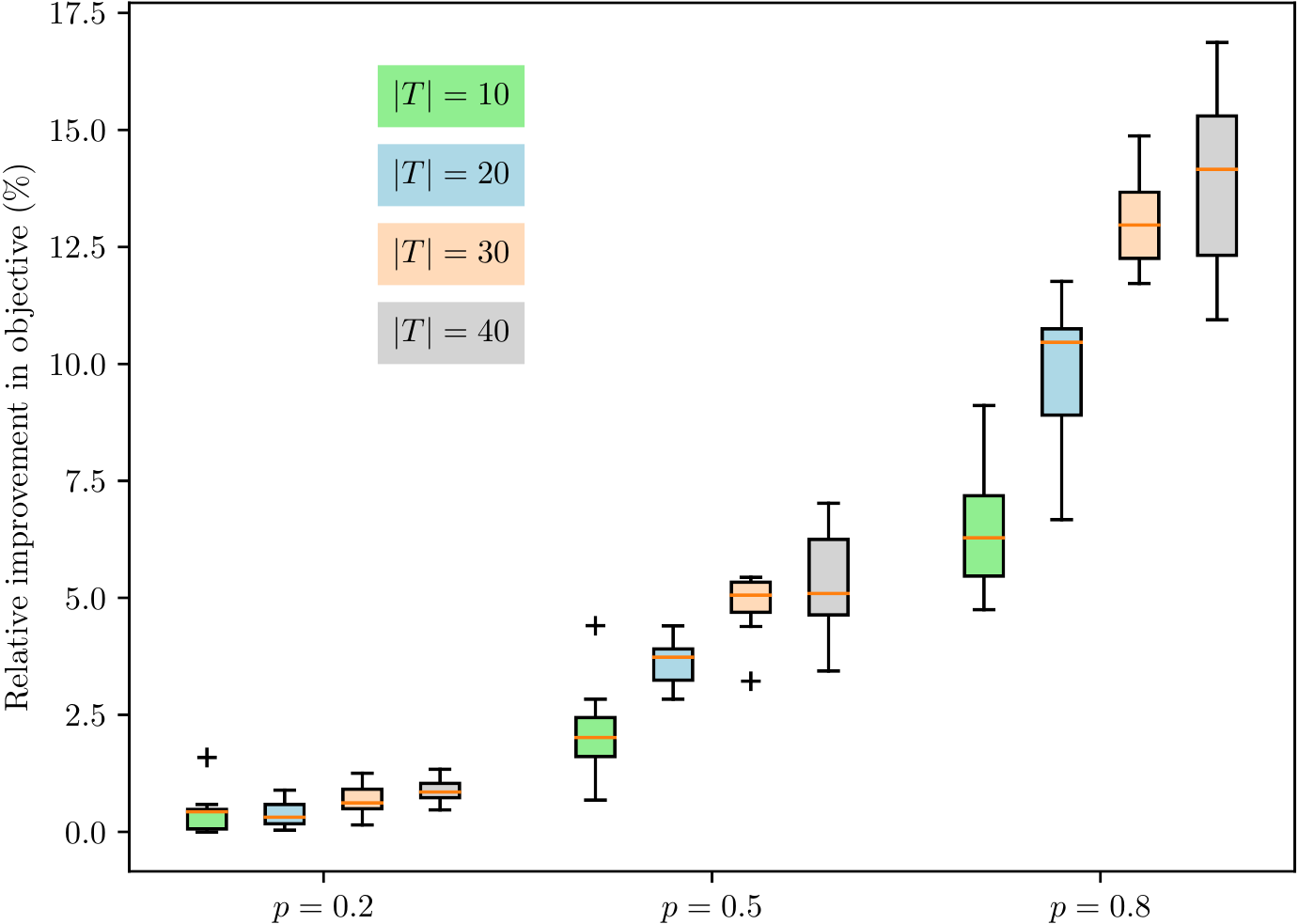}
    \caption{Value of Stochastic Solution (in percentage).}
    \label{fig:vss}
\end{figure}

\review{\paragraph{Comparison with deterministic equivalent} This set of results is focused on comparing the solution of the two-stage stochastic program provided by the PH algorithm against its the Expected Value Problem (EVP) \cite{Birge2011}. The expected value problem is solved directly using CPLEX and the value of stochastic solution or VSS is the difference between the objective value of the full two-stage stochastic program and its expected value problem. It is known that this VSS is non-negative \cite{Birge2011} and greater its value the greater the value in solving the stochastic version as opposed to the expected value problem. To this end, we fix the value of $m$ and $R$ to be $5$ units and vary the probability of generating scenarios within the set $\{0.2, 0.5, 0.8\}$. The number of generated scenarios is fixed to $100$ for all the runs. The Fig. \ref{fig:vss} show the value of stochastic solution as a relative percentage with respect to the objective value of the EVP. As observed from the Fig. \ref{fig:vss}, there is a substantial quantitative value in solving the two-stage stochastic problem as against an EVP when the value of the probabilities used for generating the scenarios are greater than $0.5$. }

\paragraph{Effect of increasing the number of scenarios} Here, we present the first set of results to demonstrate the scalability of the PH algorithm for the SVDGP with increasing number of scenarios. To that end, for this experiment, the instances with value of $n \in \{10, 15, 20, 25, 30, 35, 40\}$, $m = 5$, $R = 5$ are chosen. For each of these instances, the number of scenarios are varied in the set $|\Omega| \in \{10, 25, 50, 100, 200\}$. The results for this set of experiments is shown in Tables \ref{tab:scen-scal} and \ref{tab:scen-scal-iter}. It is observed from the tables that the PH algorithm is effectively able to solve all the instances with up to $100$ scenarios within two hours of computation time. 

\begin{table}[htbp]
    \centering
    {\onehalfspacing
    \begin{tabular}{cccccccc}
        \toprule
         $|\Omega|$ & $n = 10$ & $n = 15$ & $n = 20$ & $n = 25$ & $n = 30$ & $n=35$ & $n = 40$ \\
         \midrule
         \begin{filecontents*}{results/scenario-scalability.csv}
10,0.24,1.05,4.06,281.30,60.40,13.99,351.94
25,0.59,2.71,12.99,564.80,134.28,47.31,1084.65
50,0.95,4.85,21.86,1766.95,262.78,95.25,3121.86
100,2.02,9.94,46.57,5025.31,598.94,191.47,7256.49
200,4.28,19.07,57.45,9304.92,1133.70,336.65,16315.99
\end{filecontents*}
\csvreader[no head,late after line=\\]{results/scenario-scalability.csv}{1=\m,2=\one,3=\two,4=\three,5=\four,6=\five,7=\six,8=\seven}{\m & \one & \two & \three & \four & \five & \six & \seven}
         \bottomrule
    \end{tabular}
    \caption{Computation time (in seconds) taken by the PH algorithm for increasing number of targets and number of scenarios.}
    \label{tab:scen-scal}
    }
\end{table}

\begin{table}[htbp]
    \centering
    {\onehalfspacing
    \begin{tabular}{cccccccc}
        \toprule
         $|\Omega|$ & $n = 10$ & $n = 15$ & $n = 20$ & $n = 25$ & $n = 30$ & $n=35$ & $n = 40$ \\
         \midrule
         \begin{filecontents*}{results/scenario-scalability-iterations.csv}
10,1,1,1,16,9,1,14
25,1,1,1,14,10,1,15
50,1,1,2,26,11,1,22
100,1,1,1,34,12,1,23
200,1,1,1,30,12,1,25
\end{filecontents*}
\csvreader[no head,late after line=\\]{results/scenario-scalability-iterations.csv}{1=\m,2=\one,3=\two,4=\three,5=\four,6=\five,7=\six,8=\seven}{\m & \one & \two & \three & \four & \five & \six & \seven}
         \bottomrule
    \end{tabular}
    \caption{Number of iterations taken to obtain convergence of the PH algorithm for increasing number of targets and number of scenarios.}
    \label{tab:scen-scal-iter}
    }
\end{table}

\begin{figure}[ht!]
\begin{subfigure}{\textwidth}
  \centering
  \includegraphics[scale=0.6]{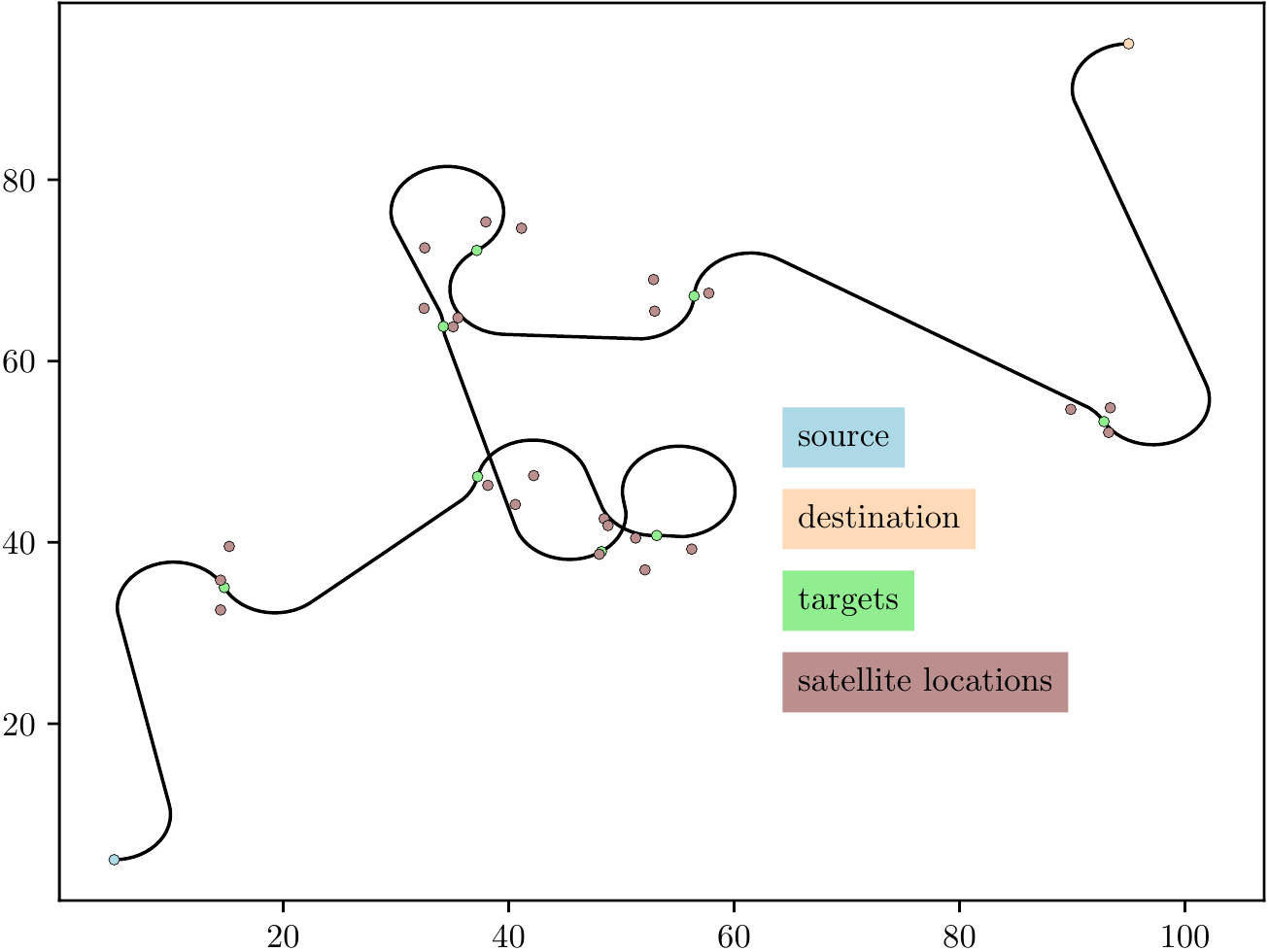}  
  \caption{First stage solution for a $10$-target instance.}
  \label{fig:sub-first}
\end{subfigure}\\
\begin{subfigure}{\textwidth}
  \centering
  \includegraphics[scale=0.6]{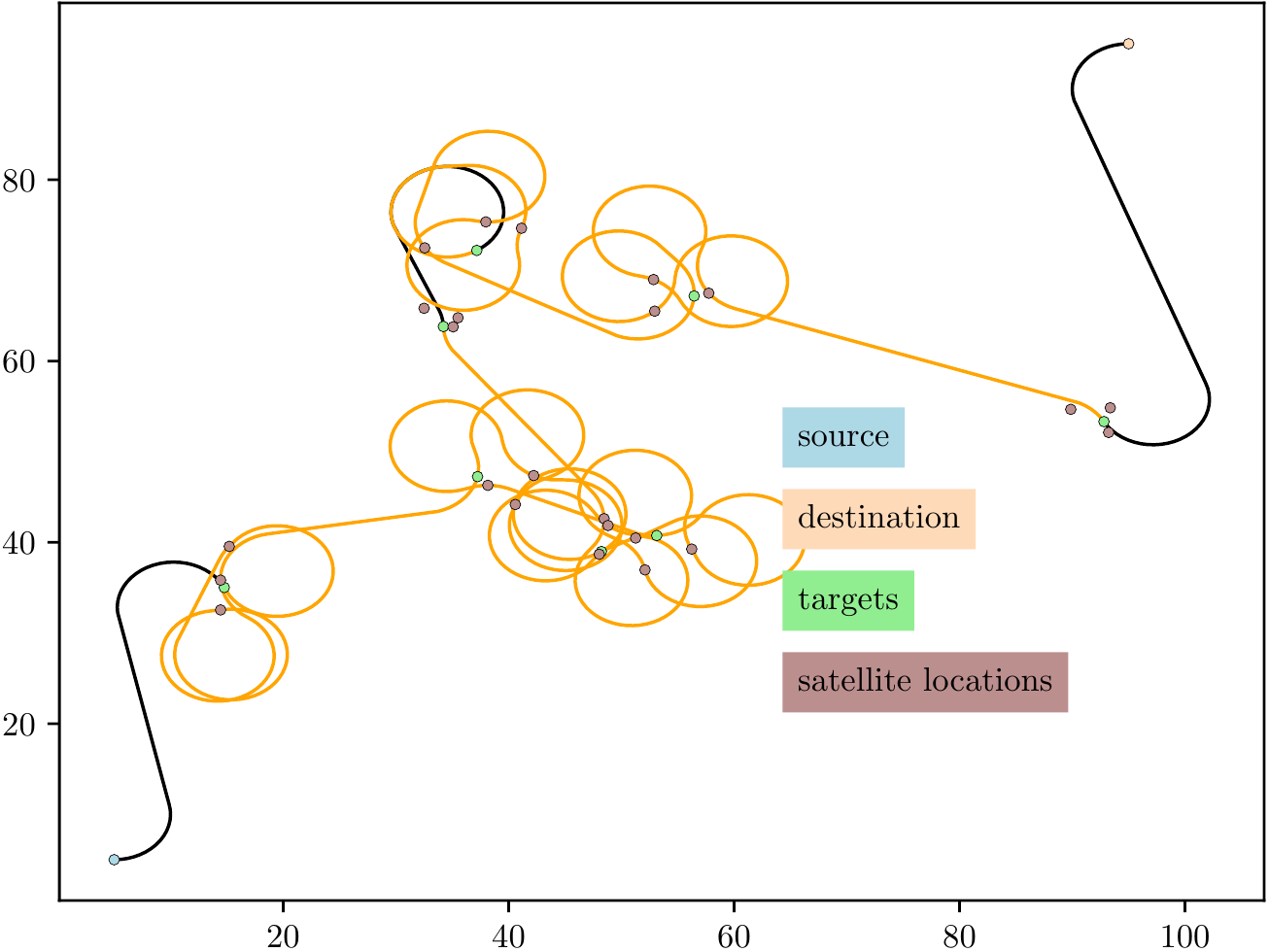} 
  \caption{The solution to the two-stage problem for one particular realization of the uncertainty. The black lines are the path that is same as the first stage solution and the orange routes correspond to the additional visits to the supplemental locations corresponding to the targets where the fidelity of information collected was not sufficient. }
  \label{fig:sub-second}
\end{subfigure}
\caption{The solution provided by the PH algorithm for a $10$-target instance. }
\label{fig:illustration}
\end{figure}

\begin{figure}[H]
    \centering
    \includegraphics[scale=0.8]{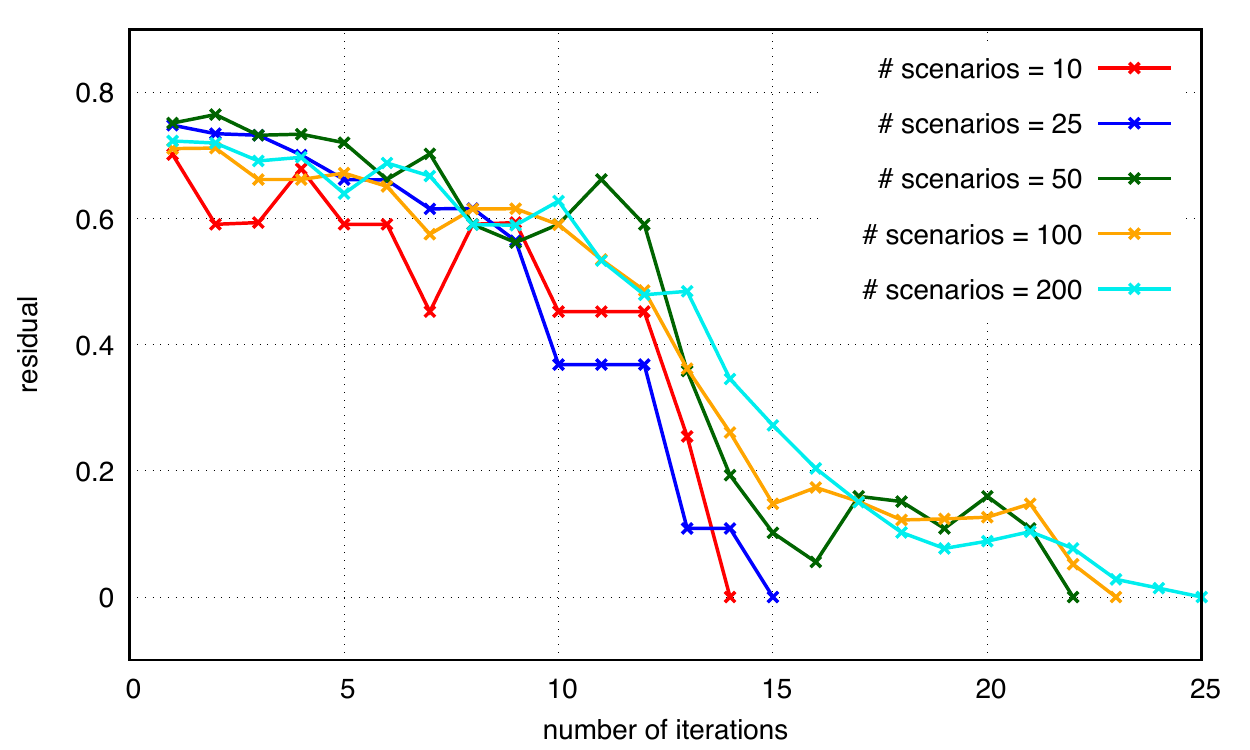}
    \caption{Residual value progress for the PH algorithm. The value of $(n, m, R)$ for the instance for this set of runs is given by $(40, 5, 5)$, respectively.}
    \label{fig:residual}
\end{figure}

The Fig. \ref{fig:residual} shows the feasible first stage (Fig. \ref{fig:sub-first}) and a second stage  solution (Fig. \ref{fig:sub-second}) for one particular realization of the random variable for a $10$-target instance. Notice that the paths taken by the vehicle are not straight lines because we assume that the vehicle is a fixed-wing aircraft with a minimum turn-radius of $5$ units and follows the Dubins path \cite{Dubins1957}. 

\paragraph{Effect of increasing the number of supplemental locations per target} The second set of results are aimed at demonstrating the scalability of the algorithm with increasing number of supplemental locations per target. To that end, for this experiment, the instances with value of $n \in \{10, 15, 20, 25, 30, 35, 40\}$, $m = \{3, 5, 7, 9\}$, $R = 5$, and $|\Omega| = 100$ (generated using a probability of $0.5$) are chosen and the results are reported in Tables \ref{tab:run-times2} and \ref{tab:iterations2}. Though the expectation is for the computation time and the number of iterations to increase with increasing value of $m$, it is not really the case for this problem, as observed from Tables \ref{tab:run-times2} and \ref{tab:iterations2}. Furthermore, we also observe that for the $40$-target instances, the computation time is much greater than two hours. \review{For typical patrolling applications, the number of target locations will be around 100 targets. For such patrolling missions, the use of multiple UAVs is preferred. In this case, a simple clustering heuristic can be used to partition the targets into manageable cluster sizes and the algorithm can be used to compute a single vehicle solution for each cluster. This addresses the scalability aspect of the algorithm. Also, since the algorithm is offline, it is typically run a few hours before the start of the mission to obtain the plans. Hence, this justifies the computation time being in the order of hours even for instances with around $40$ targets. }

\begin{table}[htbp]
    \centering
    {\onehalfspacing
    \begin{tabular}{cccccccc}
        \toprule
         $m$ & $n = 10$ & $n = 15$ & $n = 20$ & $n = 25$ & $n = 30$ & $n=35$ & $n = 40$ \\
         \midrule
         \begin{filecontents*}{results/run-times.csv}
3,1.68,6.34,21.90,185.49,92.15,66.25,2777.89
5,2.02,9.94,46.57,5025.31,598.94,191.47,7256.49
7,5.82,33.55,63.68,162.63,1605.04,6082.88,9404.69
9,22.38,57.53,870.38,1422.58,412.64,7152.56,2325.88
\end{filecontents*}
\csvreader[no head,late after line=\\]{results/run-times.csv}{1=\m,2=\one,3=\two,4=\three,5=\four,6=\five,7=\six,8=\seven}{\m & \one & \two & \three & \four & \five & \six & \seven}
         \bottomrule
    \end{tabular}
    \caption{Computation time (in seconds) taken by the PH algorithm for increasing number of targets and number of supplemental locations per target.}
    \label{tab:run-times2}
    }
\end{table}

\begin{table}[htbp]
    \centering
    {\onehalfspacing
    \begin{tabular}{cccccccc}
        \toprule
         $m$ & $n = 10$ & $n = 15$ & $n = 20$ & $n = 25$ & $n = 30$ & $n=35$ & $n = 40$ \\
         \midrule
         \begin{filecontents*}{results/iterations.csv}
3,1,1,1,11,1,1,20
5,1,1,2,34,12,1,23
7,1,1,1,1,7,19,20
9,1,1,6,6,1,8,1
\end{filecontents*}
\csvreader[no head,late after line=\\]{results/iterations.csv}{1=\m,2=\one,3=\two,4=\three,5=\four,6=\five,7=\six,8=\seven}{\m & \one & \two & \three & \four & \five & \six & \seven}
         \bottomrule
    \end{tabular}
    \caption{Number of iterations taken to obtain convergence of the PH algorithm for increasing number of targets and number of supplemental locations per target.}
    \label{tab:iterations2}
    }
\end{table}

\paragraph{Effect of increasing the pre-specifed radius for the supplemental locations} This set of results presents the variation in the computation time and the number of iterations for the PH algorithm to converge when the pre-specified radius in which the supplemental locations are present is varied in the set $\{5, 10\}$. The results present the change in computation time and number of iterations for $n = \{10, 15, 20, 25, 30, 35, 40\}$ and when $m = 5$, $|\Omega| = 100$ (generated using a probability of $0.5$). The Tables \ref{tab:run-times-R} and \ref{tab:iterations-R} present the computation times and the number of iterations for the above instances, respectively. Again, no clear trend exists between the value of $R$ and the computation time or the number of iterations taken for the PH algorithm to converge; nevertheless, the trivial trend that can be observed from the two tables is that greater the number of iterations the greater the computation time.

\begin{table}[htbp]
    \centering
    {\onehalfspacing
    \begin{tabular}{cccccccc}
        \toprule
         $R$ & $n = 10$ & $n = 15$ & $n = 20$ & $n = 25$ & $n = 30$ & $n=35$ & $n = 40$ \\
         \midrule
         \begin{filecontents*}{results/run-times-R.csv}
5,2.023,9.943,46.569,5025.308,598.943,191.472,7256.488
10,45.824,762.224,23.985,98.602,736.891,1284.960,2022.578
\end{filecontents*}
\csvreader[no head,late after line=\\]{results/run-times-R.csv}{1=\m,2=\one,3=\two,4=\three,5=\four,6=\five,7=\six,8=\seven}{\m & \one & \two & \three & \four & \five & \six & \seven}
         \bottomrule
    \end{tabular}
    \caption{Computation time (in seconds) taken by the PH algorithm for values of $R$.}
    \label{tab:run-times-R}
    }
\end{table}

\begin{table}[htbp]
    \centering
    {\onehalfspacing
    \begin{tabular}{cccccccc}
        \toprule
         $R$ & $n = 10$ & $n = 15$ & $n = 20$ & $n = 25$ & $n = 30$ & $n=35$ & $n = 40$ \\
         \midrule
         \begin{filecontents*}{results/iterations-R.csv}
5,1,1,2,34,12,1,23
10,13,37,1,1,10,12,7
\end{filecontents*}
\csvreader[no head,late after line=\\]{results/iterations-R.csv}{1=\m,2=\one,3=\two,4=\three,5=\four,6=\five,7=\six,8=\seven}{\m & \one & \two & \three & \four & \five & \six & \seven}
         \bottomrule
    \end{tabular}
    \caption{Number of iterations taken by the PH algorithm for varying values of $R$.}
    \label{tab:iterations-R}
    }
\end{table}

\paragraph{Effect of changing the probabilities in the scenario generation} This set of results is aimed at examining the effect of changing the probabilities that the information collected at any target is not of sufficient fidelity. We choose a specific instance with $n = 20$, $m = 5$, and $R=5$ units. There is no specific reason for this choice and the results for all the other instances followed the same trend.  To that end, the probabilities that we choose to generate $100$ scenarios for this particular instance is given by $\{0.0, 0.2, 0.5, 0.8, 1.0\}$. We remark that when the probability is $0.0$, it basically means that the information collected at all the targets is of sufficient fidelity and this reduces the SVDGP to computing a traveling salesman tour through the set of targets. On the other hand, if all the probabilities take a value $1.0$, then the information collected at every target is not of sufficient fidelity and the SVDGP reduces to a TSP the set of targets and their supplemental locations with additional sequencing constraints i.e., the supplemental locations have to be visited immediately after their respective target visits. The Table \ref{tab:run-time-prob} presents both the computation time and the number of iterations for this experiment. The computation time for the case when the probability value is $0.0$ is the least as solving a TSP with $20$ targets is substantially fast. On the other end of the spectrum, the computation time for case with a probability value of $1.0$ is the maximum as it reduces to the problem of solving a TSP with $120$ vertices and additional sequencing constraints. Also, for these two cases, the recourse action is obtained directly by solving the TSP or the TSP with additional sequencing constraints and hence, the number of iterations for the PH algorithm to converge for these two cases will always be equal to one.

\begin{table}[htbp]
    \centering
    {\onehalfspacing
    \begin{tabular}{ccc}
        \toprule
         $p$ & time (seconds) & \# iter \\
         \midrule
         \begin{filecontents*}{results/prob-run-times.csv}
0,4.88,1
0.2,39.62,5
0.5,45.79,2
0.8,89.38,1
1,269.26,1
\end{filecontents*}
\csvreader[no head,late after line=\\]{results/prob-run-times.csv}{1=\p,2=\t,3=\i}{\p & \t & \i}
         \bottomrule
    \end{tabular}
    \caption{Computation time and number of iterations for different values of probabilities. Here $p$ represents the probability that the information collected at the target is not of sufficient fidelity for any target, and \# iter is the number of iterations taken by the PH algorithm to converge.}
    \label{tab:run-time-prob}
    }
\end{table}

\paragraph{PH algorithm residual statistics} This set of results shows the changes in the residual of the PH algorithm (i.e., $\epsilon^k$ in Algorithm \ref{algo:ph})  for varying number of scenarios. We remark that this residual function is not a strictly decreasing function. For this experiment, the instance with $n=40$, $m=5$, and $R=5$ was chosen. As far as the scenarios are concerned, all of them were created using a probability value of $0.5$. The Fig. \ref{fig:residual} shows the residual values for varying number of scenarios. All the results presented thus far illustrate the effectiveness of the PH algorithm in computing heuristic solutions to the SVDGP with uncertainty in the fidelity of information collected. 

\section{Conclusion} \label{sec:conclusion} 
\review{This article presents the first formulation and solution approach for a UAV path planning problem in the context of data gathering missions.  This problem can be extended to a multi-stage setting using the same framework. This is the first attempt to solve such a sequential stochastic programming problem in this context and there is scope to work on this further to develop and refine distributed computing algorithms to solve it for a multi-stage setting while reducing the computational times. Furthermore, extension of this problem to a multi-vehicle setting entails addressing an additional complexity of partitioning the targets. This partitioning problem can be solved heuristically using clustering algorithms or one can resort to exact resolution of the partitions by embedding constraints in to the model itself; this presents an interesting avenue for future work.} 

\section*{Acknowledgements}
Kaarthik Sundar would like to acknowledge the funding provided by  LANL’s  Directed  Research  and  Development  (LDRD)  project: ``20200603ECR: Distributed Algorithms for Large-Scale Ordinary Differential/Partial Differential Equation (ODE/PDE) Constrained Optimization Problems on Graphs''. This work was carried out under the U.S. DOE Contract No. DE-AC52-06NA25396.

\Urlmuskip=0mu plus 1mu\relax
\bibliographystyle{elsarticle-num-names}
\bibliography{sample.bib}

\begin{thebibliography}{40}
\expandafter\ifx\csname natexlab\endcsname\relax\def\natexlab#1{#1}\fi
\providecommand{\url}[1]{\texttt{#1}}
\providecommand{\href}[2]{#2}
\providecommand{\path}[1]{#1}
\providecommand{\DOIprefix}{doi:}
\providecommand{\ArXivprefix}{arXiv:}
\providecommand{\URLprefix}{URL: }
\providecommand{\Pubmedprefix}{pmid:}
\providecommand{\doi}[1]{\href{http://dx.doi.org/#1}{\path{#1}}}
\providecommand{\Pubmed}[1]{\href{pmid:#1}{\path{#1}}}
\providecommand{\bibinfo}[2]{#2}
\ifx\xfnm\relax \def\xfnm[#1]{\unskip,\space#1}\fi
\bibitem[{Curry et~al.(2004)Curry, Maslanik, Holland, and Pinto}]{Curry2004}
\bibinfo{author}{J.~Curry}, \bibinfo{author}{J.~Maslanik},
  \bibinfo{author}{G.~Holland}, \bibinfo{author}{J.~Pinto},
\newblock \bibinfo{title}{Applications of {A}erosondes in the {A}rctic},
\newblock \bibinfo{journal}{Bulletin of the American Meteorological Society}
  \bibinfo{volume}{85} (\bibinfo{year}{2004}) \bibinfo{pages}{1855--1862}.
\bibitem[{Corrigan et~al.(2007)Corrigan, Roberts, Ramana, Kim, and
  Ramanathan}]{Corrigan2007}
\bibinfo{author}{C.~Corrigan}, \bibinfo{author}{G.~Roberts},
  \bibinfo{author}{M.~Ramana}, \bibinfo{author}{D.~Kim},
  \bibinfo{author}{V.~Ramanathan},
\newblock \bibinfo{title}{Capturing vertical profiles of aerosols and black
  carbon over the {I}ndian {O}cean using autonomous unmanned aerial vehicles},
\newblock \bibinfo{journal}{Atmospheric Chemistry and Physics Discussions}
  \bibinfo{volume}{7} (\bibinfo{year}{2007}) \bibinfo{pages}{11429--11463}.
\bibitem[{Tokekar et~al.(2016)Tokekar, Vander~Hook, Mulla, and
  Isler}]{Tokekar2016}
\bibinfo{author}{P.~Tokekar}, \bibinfo{author}{J.~Vander~Hook},
  \bibinfo{author}{D.~Mulla}, \bibinfo{author}{V.~Isler},
\newblock \bibinfo{title}{Sensor planning for a symbiotic {UAV} and {UGV}
  system for precision agriculture},
\newblock \bibinfo{journal}{IEEE Transactions on Robotics} \bibinfo{volume}{32}
  (\bibinfo{year}{2016}) \bibinfo{pages}{1498--1511}.
\bibitem[{Shanahan et~al.(2001)Shanahan, Schepers, Francis, Varvel, Wilhelm,
  Tringe, Schlemmer, and Major}]{Shanahan2001}
\bibinfo{author}{J.~F. Shanahan}, \bibinfo{author}{J.~S. Schepers},
  \bibinfo{author}{D.~D. Francis}, \bibinfo{author}{G.~E. Varvel},
  \bibinfo{author}{W.~W. Wilhelm}, \bibinfo{author}{J.~M. Tringe},
  \bibinfo{author}{M.~R. Schlemmer}, \bibinfo{author}{D.~J. Major},
\newblock \bibinfo{title}{Use of remote-sensing imagery to estimate corn grain
  yield},
\newblock \bibinfo{journal}{Agronomy Journal} \bibinfo{volume}{93}
  (\bibinfo{year}{2001}) \bibinfo{pages}{583--589}.
\bibitem[{Casbeer et~al.(2005)Casbeer, Beard, McLain, Li, and
  Mehra}]{Casbeer2005}
\bibinfo{author}{D.~W. Casbeer}, \bibinfo{author}{R.~W. Beard},
  \bibinfo{author}{T.~W. McLain}, \bibinfo{author}{S.-M. Li},
  \bibinfo{author}{R.~K. Mehra},
\newblock \bibinfo{title}{Forest fire monitoring with multiple small {UAV}s},
\newblock in: \bibinfo{booktitle}{Proceedings of the 2005, American Control
  Conference, 2005.}, \bibinfo{organization}{IEEE}, \bibinfo{year}{2005}, pp.
  \bibinfo{pages}{3530--3535}.
\bibitem[{Ferreira et~al.(2009)Ferreira, Almeida, Martins, Almeida, Dias, Dias,
  and Silva}]{Ferreira2009}
\bibinfo{author}{H.~Ferreira}, \bibinfo{author}{C.~Almeida},
  \bibinfo{author}{A.~Martins}, \bibinfo{author}{J.~Almeida},
  \bibinfo{author}{N.~Dias}, \bibinfo{author}{A.~Dias},
  \bibinfo{author}{E.~Silva},
\newblock \bibinfo{title}{Autonomous bathymetry for risk assessment with {ROAZ}
  robotic surface vehicle},
\newblock in: \bibinfo{booktitle}{Oceans 2009-Europe},
  \bibinfo{organization}{IEEE}, \bibinfo{year}{2009}, pp.
  \bibinfo{pages}{1--6}.
\bibitem[{Zajkowski et~al.(2006)Zajkowski, Dunagan, and Eilers}]{Zajkowski2006}
\bibinfo{author}{T.~Zajkowski}, \bibinfo{author}{S.~Dunagan},
  \bibinfo{author}{J.~Eilers}, \bibinfo{title}{Small {UAS} communications
  mission}, \bibinfo{year}{2006}. \URLprefix
  \url{https://www.fs.fed.us/eng/rsac/RS2006/presentations/zajkowski2.pdf}.
\bibitem[{Manyam et~al.(2019)Manyam, Sundar, and Casbeer}]{Manyam2019}
\bibinfo{author}{S.~G. Manyam}, \bibinfo{author}{K.~Sundar},
  \bibinfo{author}{D.~W. Casbeer},
\newblock \bibinfo{title}{Cooperative routing for an air-ground vehicle
  team--exact algorithm, transformation method, and heuristics},
\newblock \bibinfo{journal}{IEEE Transactions on Automation Science and
  Engineering}  (\bibinfo{year}{2019}).
\bibitem[{Sundar et~al.(2017)Sundar, Venkatachalam, and Manyam}]{Sundar2017}
\bibinfo{author}{K.~Sundar}, \bibinfo{author}{S.~Venkatachalam},
  \bibinfo{author}{S.~G. Manyam},
\newblock \bibinfo{title}{Path planning for multiple heterogeneous unmanned
  vehicles with uncertain service times},
\newblock in: \bibinfo{booktitle}{Unmanned Aircraft Systems (ICUAS), 2017
  International Conference on}, \bibinfo{organization}{IEEE},
  \bibinfo{year}{2017}, pp. \bibinfo{pages}{480--487}.
\bibitem[{Shima and Rasmussen(2009)}]{Shima2009}
\bibinfo{author}{T.~Shima}, \bibinfo{author}{S.~Rasmussen}, \bibinfo{title}{UAV
  cooperative decision and control: challenges and practical approaches},
  \bibinfo{publisher}{SIAM}, \bibinfo{year}{2009}.
\bibitem[{Liu et~al.(2019)Liu, Liu, Shi, Wu, and Chen}]{Liu2019}
\bibinfo{author}{Y.~Liu}, \bibinfo{author}{Z.~Liu}, \bibinfo{author}{J.~Shi},
  \bibinfo{author}{G.~Wu}, \bibinfo{author}{C.~Chen},
\newblock \bibinfo{title}{{Optimization of Base Location and Patrol Routes for
  Unmanned Aerial Vehicles in Border Intelligence, Surveillance, and
  Reconnaissance}},
\newblock \bibinfo{journal}{Journal of Advanced Transportation}
  \bibinfo{volume}{2019} (\bibinfo{year}{2019}).
\bibitem[{Adams and Friedland(2011)}]{Adams2011}
\bibinfo{author}{S.~M. Adams}, \bibinfo{author}{C.~J. Friedland},
\newblock \bibinfo{title}{{A survey of unmanned aerial vehicle (UAV) usage for
  imagery collection in disaster research and management}},
\newblock in: \bibinfo{booktitle}{9th International Workshop on Remote Sensing
  for Disaster Response}, \bibinfo{year}{2011}, p.~\bibinfo{pages}{8}.
\bibitem[{Rajan et~al.(2019)Rajan, Sundar, and Gautam}]{rajan2019routing}
\bibinfo{author}{S.~Rajan}, \bibinfo{author}{K.~Sundar},
  \bibinfo{author}{N.~Gautam},
\newblock \bibinfo{title}{Routing problems for reconnaissance patrolling
  missions},
\newblock in: \bibinfo{booktitle}{2019 International Conference on Unmanned
  Aircraft Systems (ICUAS)}, \bibinfo{organization}{IEEE},
  \bibinfo{year}{2019}, pp. \bibinfo{pages}{213--220}.
\bibitem[{Balasubramanian and Grossmann(2004)}]{Balasubramanian2004}
\bibinfo{author}{J.~Balasubramanian}, \bibinfo{author}{I.~Grossmann},
\newblock \bibinfo{title}{Approximation to multistage stochastic optimization
  in multiperiod batch plant scheduling under demand uncertainty},
\newblock \bibinfo{journal}{Industrial \& engineering chemistry research}
  \bibinfo{volume}{43} (\bibinfo{year}{2004}) \bibinfo{pages}{3695--3713}.
\bibitem[{Sundar and Rathinam(2013)}]{Sundar2013}
\bibinfo{author}{K.~Sundar}, \bibinfo{author}{S.~Rathinam},
\newblock \bibinfo{title}{Algorithms for routing an unmanned aerial vehicle in
  the presence of refueling depots},
\newblock \bibinfo{journal}{IEEE Transactions on Automation Science and
  Engineering} \bibinfo{volume}{11} (\bibinfo{year}{2013})
  \bibinfo{pages}{287--294}.
\bibitem[{Sundar and Rathinam(2017)}]{Sundar2017JINT}
\bibinfo{author}{K.~Sundar}, \bibinfo{author}{S.~Rathinam},
\newblock \bibinfo{title}{Algorithms for heterogeneous, multiple depot,
  multiple unmanned vehicle path planning problems},
\newblock \bibinfo{journal}{Journal of Intelligent \& Robotic Systems}
  \bibinfo{volume}{88} (\bibinfo{year}{2017}) \bibinfo{pages}{513--526}.
\bibitem[{Levy et~al.(2014)Levy, Sundar, and Rathinam}]{Levy2014}
\bibinfo{author}{D.~Levy}, \bibinfo{author}{K.~Sundar},
  \bibinfo{author}{S.~Rathinam},
\newblock \bibinfo{title}{Heuristics for routing heterogeneous unmanned
  vehicles with fuel constraints},
\newblock \bibinfo{journal}{Mathematical Problems in Engineering}
  \bibinfo{volume}{2014} (\bibinfo{year}{2014}).
\bibitem[{Otto et~al.(2018)Otto, Agatz, Campbell, Golden, and Pesch}]{Otto2018}
\bibinfo{author}{A.~Otto}, \bibinfo{author}{N.~Agatz},
  \bibinfo{author}{J.~Campbell}, \bibinfo{author}{B.~Golden},
  \bibinfo{author}{E.~Pesch},
\newblock \bibinfo{title}{Optimization approaches for civil applications of
  unmanned aerial vehicles (uavs) or aerial drones: A survey},
\newblock \bibinfo{journal}{Networks} \bibinfo{volume}{72}
  (\bibinfo{year}{2018}) \bibinfo{pages}{411--458}.
\bibitem[{Gendreau et~al.(1996)Gendreau, Laporte, and
  S{\'e}guin}]{Gendreau1996}
\bibinfo{author}{M.~Gendreau}, \bibinfo{author}{G.~Laporte},
  \bibinfo{author}{R.~S{\'e}guin},
\newblock \bibinfo{title}{Stochastic vehicle routing},
\newblock \bibinfo{journal}{European Journal of Operational Research}
  \bibinfo{volume}{88} (\bibinfo{year}{1996}) \bibinfo{pages}{3--12}.
\bibitem[{Jaillet(1988)}]{Jaillet1988}
\bibinfo{author}{P.~Jaillet},
\newblock \bibinfo{title}{A priori solution of a traveling salesman problem in
  which a random subset of the customers are visited},
\newblock \bibinfo{journal}{Operations research} \bibinfo{volume}{36}
  (\bibinfo{year}{1988}) \bibinfo{pages}{929--936}.
\bibitem[{Carraway et~al.(1989)Carraway, Morin, and Moskowitz}]{Carraway1989}
\bibinfo{author}{R.~L. Carraway}, \bibinfo{author}{T.~L. Morin},
  \bibinfo{author}{H.~Moskowitz},
\newblock \bibinfo{title}{Generalized dynamic programming for stochastic
  combinatorial optimization},
\newblock \bibinfo{journal}{Operations Research} \bibinfo{volume}{37}
  (\bibinfo{year}{1989}) \bibinfo{pages}{819--829}.
\bibitem[{Tillman(1969)}]{Tillman1969}
\bibinfo{author}{F.~A. Tillman},
\newblock \bibinfo{title}{The multiple terminal delivery problem with
  probabilistic demands},
\newblock \bibinfo{journal}{Transportation Science} \bibinfo{volume}{3}
  (\bibinfo{year}{1969}) \bibinfo{pages}{192--204}.
\bibitem[{Dror et~al.(1989)Dror, Laporte, and Trudeau}]{Dror1989}
\bibinfo{author}{M.~Dror}, \bibinfo{author}{G.~Laporte},
  \bibinfo{author}{P.~Trudeau},
\newblock \bibinfo{title}{Vehicle routing with stochastic demands: Properties
  and solution frameworks},
\newblock \bibinfo{journal}{Transportation science} \bibinfo{volume}{23}
  (\bibinfo{year}{1989}) \bibinfo{pages}{166--176}.
\bibitem[{Bertsimas(1988)}]{Bertsimas1988}
\bibinfo{author}{D.~Bertsimas}, \bibinfo{title}{Probabilistic combinatorial
  optimization problems}, Ph.D. thesis, Massachusetts Institute of Technology,
  \bibinfo{year}{1988}.
\bibitem[{Jaillet and Odoni(1988)}]{Jaillet1988a}
\bibinfo{author}{P.~Jaillet}, \bibinfo{author}{A.~Odoni},
  \bibinfo{title}{Vehicle routing: Methods and studies, chapter the
  probabilistic vehicle routing problem}, \bibinfo{year}{1988}.
\bibitem[{Jezequel(1985)}]{Jezequel1985}
\bibinfo{author}{A.~Jezequel}, \bibinfo{title}{Probabilistic vehicle routing
  problems}, Ph.D. thesis, Massachusetts institute of technology,
  \bibinfo{year}{1985}.
\bibitem[{Shavarani et~al.(2018)Shavarani, Nejad, Rismanchian, and
  Izbirak}]{Shavarani2018}
\bibinfo{author}{S.~M. Shavarani}, \bibinfo{author}{M.~G. Nejad},
  \bibinfo{author}{F.~Rismanchian}, \bibinfo{author}{G.~Izbirak},
\newblock \bibinfo{title}{Application of hierarchical facility location problem
  for optimization of a drone delivery system: a case study of amazon prime air
  in the city of san francisco},
\newblock \bibinfo{journal}{The International Journal of Advanced Manufacturing
  Technology} \bibinfo{volume}{95} (\bibinfo{year}{2018})
  \bibinfo{pages}{3141--3153}.
\bibitem[{Torabbeigi et~al.(2018)Torabbeigi, Lim, and Kim}]{Torabbeigi2018}
\bibinfo{author}{M.~Torabbeigi}, \bibinfo{author}{G.~J. Lim},
  \bibinfo{author}{S.~J. Kim},
\newblock \bibinfo{title}{Drone delivery schedule optimization considering the
  reliability of drones},
\newblock in: \bibinfo{booktitle}{2018 International Conference on Unmanned
  Aircraft Systems (ICUAS)}, \bibinfo{organization}{IEEE},
  \bibinfo{year}{2018}, pp. \bibinfo{pages}{1048--1053}.
\bibitem[{Venkatachalam et~al.(2018)Venkatachalam, Sundar, and
  Rathinam}]{Venkatachalam2018}
\bibinfo{author}{S.~Venkatachalam}, \bibinfo{author}{K.~Sundar},
  \bibinfo{author}{S.~Rathinam},
\newblock \bibinfo{title}{A two-stage approach for routing multiple unmanned
  aerial vehicles with stochastic fuel consumption},
\newblock \bibinfo{journal}{Sensors} \bibinfo{volume}{18}
  (\bibinfo{year}{2018}) \bibinfo{pages}{3756}.
\bibitem[{Liu et~al.(2019)Liu, Liu, Zhu, and Zheng}]{Liu2019Stochastic}
\bibinfo{author}{M.~Liu}, \bibinfo{author}{X.~Liu}, \bibinfo{author}{M.~Zhu},
  \bibinfo{author}{F.~Zheng},
\newblock \bibinfo{title}{Stochastic drone fleet deployment and planning
  problem considering multiple-type delivery service},
\newblock \bibinfo{journal}{Sustainability} \bibinfo{volume}{11}
  (\bibinfo{year}{2019}) \bibinfo{pages}{3871}.
\bibitem[{Watson and Woodruff(2011)}]{Watson2011}
\bibinfo{author}{J.-P. Watson}, \bibinfo{author}{D.~L. Woodruff},
\newblock \bibinfo{title}{Progressive hedging innovations for a class of
  stochastic mixed-integer resource allocation problems},
\newblock \bibinfo{journal}{Computational Management Science}
  \bibinfo{volume}{8} (\bibinfo{year}{2011}) \bibinfo{pages}{355--370}.
\bibitem[{Rockafellar and Wets(1991)}]{Rockafellar1991}
\bibinfo{author}{R.~T. Rockafellar}, \bibinfo{author}{R.~J.-B. Wets},
\newblock \bibinfo{title}{Scenarios and policy aggregation in optimization
  under uncertainty},
\newblock \bibinfo{journal}{Mathematics of operations research}
  \bibinfo{volume}{16} (\bibinfo{year}{1991}) \bibinfo{pages}{119--147}.
\bibitem[{Shapiro(2011)}]{Shapiro2011}
\bibinfo{author}{A.~Shapiro},
\newblock \bibinfo{title}{Analysis of stochastic dual dynamic programming
  method},
\newblock \bibinfo{journal}{European Journal of Operational Research}
  \bibinfo{volume}{209} (\bibinfo{year}{2011}) \bibinfo{pages}{63--72}.
\bibitem[{Atakan and Sen(2018)}]{Atakan2018}
\bibinfo{author}{S.~Atakan}, \bibinfo{author}{S.~Sen},
\newblock \bibinfo{title}{A progressive hedging based branch-and-bound
  algorithm for mixed-integer stochastic programs},
\newblock \bibinfo{journal}{Computational Management Science}
  \bibinfo{volume}{15} (\bibinfo{year}{2018}) \bibinfo{pages}{501--540}.
\bibitem[{Fan and Liu(2010)}]{Fan2010}
\bibinfo{author}{Y.~Fan}, \bibinfo{author}{C.~Liu},
\newblock \bibinfo{title}{Solving stochastic transportation network protection
  problems using the progressive hedging-based method},
\newblock \bibinfo{journal}{Networks and Spatial Economics}
  \bibinfo{volume}{10} (\bibinfo{year}{2010}) \bibinfo{pages}{193--208}.
\bibitem[{Listes and Dekker(2005)}]{Listes2005}
\bibinfo{author}{O.~Listes}, \bibinfo{author}{R.~Dekker},
\newblock \bibinfo{title}{A scenario aggregation--based approach for
  determining a robust airline fleet composition for dynamic capacity
  allocation},
\newblock \bibinfo{journal}{Transportation Science} \bibinfo{volume}{39}
  (\bibinfo{year}{2005}) \bibinfo{pages}{367--382}.
\bibitem[{L{\o}kketangen and Woodruff(1996)}]{Lokketangen1996}
\bibinfo{author}{A.~L{\o}kketangen}, \bibinfo{author}{D.~L. Woodruff},
\newblock \bibinfo{title}{Progressive hedging and tabu search applied to mixed
  integer (0, 1) multistage stochastic programming},
\newblock \bibinfo{journal}{Journal of Heuristics} \bibinfo{volume}{2}
  (\bibinfo{year}{1996}) \bibinfo{pages}{111--128}.
\bibitem[{Wallace and Ziemba(2005)}]{Wallace2005}
\bibinfo{author}{S.~W. Wallace}, \bibinfo{author}{W.~T. Ziemba},
  \bibinfo{title}{Applications of stochastic programming},
  \bibinfo{publisher}{SIAM}, \bibinfo{year}{2005}.
\bibitem[{Dubins(1957)}]{Dubins1957}
\bibinfo{author}{L.~E. Dubins},
\newblock \bibinfo{title}{On curves of minimal length with a constraint on
  average curvature, and with prescribed initial and terminal positions and
  tangents},
\newblock \bibinfo{journal}{American Journal of mathematics}
  \bibinfo{volume}{79} (\bibinfo{year}{1957}) \bibinfo{pages}{497--516}.
\bibitem[{Birge and Louveaux(2011)}]{Birge2011}
\bibinfo{author}{J.~R. Birge}, \bibinfo{author}{F.~Louveaux},
  \bibinfo{title}{Introduction to stochastic programming},
  \bibinfo{publisher}{Springer Science \& Business Media},
  \bibinfo{year}{2011}.

\end{thebibliography}

\end{document}